\input amstex
\documentstyle{amsppt}
\documentstyle{amsppt}
\NoBlackBoxes \raggedbottom
 \topmatter


\smallskip

\NoBlackBoxes
\topmatter
\title   Surgery in codimension 3 and the Browder--Livesay invariants
\endtitle
\rightheadtext{Surgery in codimension 3}
\leftheadtext{F. Hegenbarth, Yu. V. Muranov and D. Repov\v s}
\author Friedrich Hegenbarth -- Yuri V. Muranov -- Du\v san Repov\v s
\endauthor
\thanks  This research was supported by SRA grants P1-0292-0101-04 and J1-9643-0101-07.
The second author thanks  Max-Planck Institute f\" ur Mathematik (Bonn, Germany)  for
the hospitality and the financial support.
\endthanks
\keywords Surgery assembly map, closed manifolds surgery problem, Assembly map,
inertia subgroup,
splitting problem, Browder--Livesay invariants, Browder-Livesay
groups,  normal maps, iterated  Browder--Livesay invariants, manifolds with
filtration, Browder-Quinn surgery obstruction groups, elements of
the second type of a Wall group
\endkeywords
\subjclass\nofrills{2000 {\it Mathematics Subject Classification.}
Primary 57R67, 19J25;   Secondary 55T99, 58A35, 18F25}
\endsubjclass
\abstract The inertia subgroup $I_n(\pi)$ of a surgery obstruction
group $L_n(\pi)$ is generated by elements which act trivially on
the set of homotopy triangulations $\Cal S(X)$ for some closed
topological manifold $X^{n-1}$ with $\pi_1(X)=\pi$. This group is
a subgroup of the group $C_n(\pi)$ which  consists of the elements
which  can be realized by normal maps of closed manifolds. In all
known cases  these groups coincide and the computation of them is
one of the basic problems of surgery theory.  The computation  of the
group $C_n(\pi)$ is equivalent to the computation  the image of the
assembly map $A:H_{n}(B\pi, \bold L_{\bullet})\to L_{n}(\pi)$. Every 
Browder-Livesay filtration of the manifold $X$ provides a collection of
Browder-Livesay invariants which are the forbidden invariants in the
closed manifold surgery problem.   In the present paper we
describe all possible forbidden invariants which can give a
Browder-Livesay filtration for computing the inertia subgroup. Our
approach is a natural generalization of the approach of Hambleton
and Kharshiladze. More precisely, we prove that a Browder-Livesay
filtration  of a given manifold can give the following forbidden
invariants for an element $x\in L_n(\pi_1(X))$ to belong to the subgroup
$I_n(\pi)$: the nontrivial Browder-Livesay invariants in
codimensions 0, 1, 2 and a  nontrivial class of obstructions of a restriction
of a normal map to a submanifold in codimension  3.
\endabstract
 
\endtopmatter

\document

\bigskip
\subhead 1. Introduction
\endsubhead
\bigskip

Throughout the paper we consider finitely presented groups $\pi$
equipped with an orientation  homomorphisms
$w\colon \pi\to\{\pm 1 \}$. Let $L_n(\pi)$ be  the Wall surgery
obstruction groups $L_n^s(\pi, w)$.
 As usually,
$B\pi=K(\pi,1)$ denotes the classifying space of $\pi$.
 For a manifold $X$ we suppose
that  the orientation map $w\colon \pi=\pi_1(X)\to \{\pm 1\}$
coincides with the Stiefel-Whitney character. We shall work in the
category of topological manifolds.

Any element $x\in L_{n+1}(\pi)$ is represented by a normal map of
a closed manifold with boundary. The results of \cite{17}
provide the following representation.
 Choose a closed $n$@-manifold
$X^{n}$ with $\pi_1(X^{n})=\pi$. Then there is  a normal  map
$$
(F,B)\colon (W^{n+1};\partial_0W,\partial_1W)\to (X\times I;X\times\{0\},
X\times\{1\}) \tag 1.1
$$
with the obstruction $\sigma(F,B)=x$,  where $\partial_0 W=X$,
$F|_{\partial_0W}=\operatorname{Id}$, $\partial_1W=M^{n}$, and 
$$
f=F|_M: M\to X 
$$ 
is a simple homotopy equivalence. Here
$B:\nu_W\to \nu_{X\times I}$ is a  covering $F$ bundle map of the
stable normal bundle of $W$ in Euclidean space to a bundle over
$X\times I$. In what follows we shall not mention the maps of
stable bundles if this does not lead to a confusion.

Let  $X$ be  a closed $n$@-dimensional topological manifold. An
orientation preserving  simple homotopy equivalence $f\: M^n\to  X^n$
of $n$@-manifolds
is called {\sl a  homotopical triangulation} of the manifold $X$.
Two homotopical  triangulations $f_i:M_i\to X (i=1,2)$ are
said to be equivalent if there exists a homeomorphism $h:M_1\to M_2$ fitting
in the following  homotopy commutative diagram (see \cite{14}, \cite{15}, and
\cite{17})
$$
\matrix
M_1&\overset{f_1}\to{\longrightarrow}& X. \\
\downarrow_h& \ \ \nearrow_{f_2}&\\
M_2&&\\
\endmatrix
\tag 1.2
$$
The set of equivalence classes  $\Cal S^{TOP}(X)$ fits into
the surgery exact sequence  ( see \cite{14}, \cite{15}, and \cite{17})
$$
\cdots \to[\Sigma X, G/TOP] \longrightarrow L_{n+1}(\pi)
\overset{\lambda}\to{\longrightarrow} \Cal S^{TOP}(X)
\longrightarrow [X, G/TOP] \overset{\sigma}\to{\longrightarrow}
L_n(\pi). \tag 1.3
$$
\smallskip

 For a homotopy
triangulation $g\colon M\to X$, the
action $\lambda$ in (1.3)  is defined in a similar way  as in  the
representation (1.1) that gives the  action of $x\in L_{n+1}(\pi)$
on the trivial triangulation $\operatorname{Id}\colon  X\to X$. By
definition,  $\left[\lambda(x)\right](\operatorname{Id})$ is the homotopy triangulation
$$
F|_{\partial_1W}\colon  \partial_1 W\to X\times\{1\}
$$
of the manifold $X$ \cite{17, \S 10}. Let an element $x\in
L_{n+1}(\pi)$ act trivially  on a homotopy triangulation of some
manifold $X^n$ with $\pi_1(X)=\pi$. We denote by $I_{n+1}(\pi)$
the subgroup generated by such elements.

In fact, Wall \cite{17} constructed an action of an element $x\in
L_{n+1}(\pi)$ on any homotopy triangulation
$$
f\colon (M,\partial M)\to (X, \partial X)
\tag 1.4
$$
of compact  manifolds relative the boundary. Let
$\pi_1(X^n)=\pi_1(\partial X)=\pi$. Then the results of Wall provide a
normal map of $4$-ads
$$
F\colon (W^{n+1}; \partial_0W, \partial_1W, V) \to (X\times I; X\times{0},
X\times{1}, \partial X\times I)
\tag 1.5
$$
with $\sigma(F)=x$. In (1.5) we have $\partial_0W=X\times{0}$,
$V=\partial X\times I$, $F|_{\partial_0W \cup
V}=\operatorname{Id}$, $\partial_1W=M$ with the boundary $\partial M=
\partial X$,  and the map
$$
f_1=F|_{\partial_1W}\colon \partial_1W\to X\times{1}\tag 1.6
$$
is a simple homotopy equivalence.
\smallskip

The assembly map
$$
A:H_{n}(B\pi; \bold L_{\bullet})\to L_{n}(\pi) \tag 1.7
$$
fits in the following algebraic surgery exact sequence of Ranicki \cite{14}
for the space $B\pi$
$$
\cdots \to L_{n+1}(\pi) \rightarrow \Cal S_{n+1}(B\pi) \to
H_{n}(B\pi; \bold L_{\bullet})  \overset{A}\to{\rightarrow}
L_n(\pi) \to \cdots \tag 1.8
$$
where $\Cal S_{n+1}(B\pi)$ is the structure set of the topological space $B\pi$ and
 $\bold L_{\bullet}$ is the  $\Omega$@-spectrum
that is  an 1-connected cover of the simply connected surgery
$\Omega$@-spectrum $\bold L_{\bullet}(1)$ with  $\pi_n( \bold
L_{\bullet}(1))= L_n(1) \ (n>0)$ and $\bold {L_{\bullet}}_0 \simeq G/TOP$
 \cite{14}.
For a closed  $n$@-dimensional topological manifold $X$  the
surgery  exact sequence (1.3) is isomorphic to the left part (from
the group $L_n(\pi)$) of the algebraic surgery exact sequence
(1.8) of the space $X$  with
$$
H_n (X; \bold L_{\bullet}) \cong  [X, G / TOP] \ \ \text{and} \ \
\Cal S_{n+1} (X) = \Cal{S}^{TOP} (X)
$$
(see \cite{14}, \cite{15}, and \cite{17}).
The image $C_n(\pi)\subset L_n(\pi)$ of the assembly map consists
of the surgery obstructions $\sigma(f,b)\in L_n(\pi)$ where
$(f,b)$ is a normal map of closed manifolds  with a given
orientation map $w$ \cite{17, \S 13B}. It was proved in  \cite{5}
that $I_{n}(\pi)\subset C_{n}(\pi)$.
Additionally, it follows from \cite{7} that images of these groups
coincide in the projective Novikov-Wall groups $L^p_*$. The problem of
the computation of the groups $I_n(\pi)$ and $C_n(\pi)$  is one of
the significant problems of geometrical topology (see \cite{8}).

The iterated Browder-Livesay invariants  provide a collection of
the forbidden invariants in the
closed manifold surgery problem for a group $\pi$ with a subgroup
of index 2 (see \cite{1}, \cite{3}, \cite{5}, \cite{7}, \cite{11}, and \cite{12}).
The natural way to describe  iterated Browder-Livesay invariants for the group
$L_n(\pi)$ is
the Browder-Livesay filtration  of the manifold $X$ with $\pi_1(X)=\pi$
(see \cite{1}, \cite{6},  and \cite{13}) (the definition of a Browder-Livesay filtration and iterated invariants is given in section 2
below).
  In the present paper we describe the application of this approach to
computing  the inertia subgroup.
 More precisely, we prove that a Browder-Livesay
filtration  of a given manifold can give the following forbidden
invariants for an element $x\in L_n(\pi_1(X))$ to belong to the subgroup
$I_n(\pi)$: the nontrivial Browder-Livesay invariants in
codimensions 0, 1, 2 and a  nontrivial class of obstructions of a restriction
of a normal map to a submanifold in codimension  3.

In section 2 we give necessary preliminary material and state the main theorem  and in section 3
we formulae and prove the main results. 
\bigskip

\subhead 2. Browder-Livesay filtration and iterated invariants
\endsubhead
\bigskip

 The pair
of closed manifolds $(X^n,Y^{n-1})$ is called {\sl a Browder-Livesay pair}
if  $Y$ is the one-sided
submanifold in $X$ and induced by the natural inclusion map
$\pi_1(Y)\to \pi_1(X)$ is an isomorphism (see \cite{2}, \cite{4}, \cite{5},  \cite{7}, \cite{12}, and \cite{17, \S 11}).
Let $\rho=\pi_1(X\setminus Y)$ and $i\colon \rho\to \pi$ be the natural map induced by the
inclusion.
  Let $U$ be a tubular neighborhood
 of $Y$ in $X$ with a boundary $\partial U$.
 We obtain  a pushout square
$$
F=
\left(\matrix
\pi_1(\partial U)& \to &\pi_1(X\setminus Y)\\
\downarrow & & \downarrow\\
\pi_1(Y)&\to &\pi_1(X)\\
\endmatrix\right) =
\left(\matrix
\rho& \to &\rho\\
\downarrow & & \downarrow\\
\pi^{\mp}&\to &\pi^{\pm}\\
\endmatrix\right)
\tag 2.1
$$
of fundamental groups with orientation in which horizontal maps
are isomorphisms and vertical maps are inclusions of index 2. In (2.1) the upper horizontal
 map and the vertical maps  agree with orientations. The bottom horizontal
 maps preserve the orientation on the image of the vertical maps and
 reverse orientations outside these images.
We shall denote this fact by superscript "$+$" or "$-$". We shall
omit this superscript if  the orientation follows from the
context.

By definition (see \cite{15, \S 7} and \cite{17, \S
11}), a
 simple homotopy equivalence $f: M\to X$ {\sl splits along the submanifold $Y$}
if it is homotopic to a  map $g : M\to X$ which is transversal to
$Y$ with $N = g^{-1}(Y)$, and whose restrictions
$$
g|_N : N\to Y \ \ \text{and}\ \ g|_{(M \setminus N)}:M\setminus
N\to X\setminus Y \tag 2.2
$$
are simple homotopy equivalences. The splitting
obstruction groups
$$
LN_{n-1}(\pi_1(X\setminus Y)\to \pi_1(X))= LN_{n-1}(\rho\to \pi)
$$
 for a Browder-Livesay manifold pair $(X,Y)$ were defined
(see \cite{2}, \cite{3}, \cite{5}, \cite{7}, \cite{15, \S 7},
and \cite{17}) and are called  the {\sl Browder--Livesay groups}.
 The algebraic definition  of the
$LN_*$@-groups is given in \cite{16}. These groups depend
functorially on the oriented inclusion $\rho\to \pi$
 and the dimension $n-q\bmod 4$.

These groups fit in the following braid of exact sequences
 ( see \cite{2}, \cite{7},  \cite{12}, \cite{15},  and  \cite{16})
$$
 \matrix \rightarrow & {L}_{n}(\rho) & \overset{i_*}\to{\longrightarrow} &
   L_{n}(\pi) &
\overset{\partial}\to{\rightarrow} & LN_{n-2}(\rho\to
\pi)&\rightarrow
 \cr
\ & \nearrow \  \  \ \ \ \ \  \ \searrow &\ &s\nearrow \ \ \ \ \
                           \searrow
   & \ & \nearrow \ \ \  \ \ \ \ \ \ \searrow & \ \cr
  \ & \ & LP_{n-1}(F)& \Gamma\downarrow & L_{n}(i_*) & \ & \ \cr
   \ & \searrow \ \ \ \ \ \ \ \ \ \nearrow &\ &q\searrow  \ \ \ \ \
\nearrow
 & \ & \searrow \ \   \ \ \ \  \ \ \ \nearrow & \ \cr
\rightarrow & LN_{n-1}(\rho\to \pi) &
\overset{c}\to{\longrightarrow} &
 L_{n-1}(\pi^-) &
 \overset{i^!_-}\to{\longrightarrow} & {L}_{n-1}(\rho) & \rightarrow
\endmatrix
\tag 2.3
 $$
 \smallskip
\noindent where $LP_{n-1}(F)\cong L_n(i^!_-)$ are the surgery
obstruction groups for the manifold pair $(X,Y)$ (see \cite{2},
\cite{15}, and \cite{17}),  and $L_n(i_*)$ are the relative surgery obstruction groups
for the inclusion $i$ (see \cite{2}, \cite{10}, \cite{15},  and \cite{17}).
The upper and bottom rows of Diagram (2.3) are
chain complexes, and  $\Gamma$ is an  isomorphism of the
corresponding homology groups.
Note that  the
maps $s$ and $q$ are the natural forgetful maps,  and  the map $c$
denote passing from surgery problem inside the manifold $X$ to
an abstract surgery problem \cite{17}. The map $i^!_-$ is the surgery transfer
map, and the map $\partial$ is the composition
$$
L_{n}(\pi)\overset{\lambda}\to{\to} \Cal S^{TOP}(X)\to LN_{n-2}(\rho\to\pi)
 $$
of the action of an element $x$ on the trivial triangulation of a
closed manifold $X^{n-1}$ and taking an obstruction to splitting
along the submanifold $Y^{n-2}\subset X$  on the top boundary of
the bordism as in (1.1).  For an element $x\in L_{n}(\pi)$
which  represents a homology class
$$
[x]\in \operatorname{Ker}\partial/\operatorname{Im}i_*
$$
we have a class
$$
\Gamma([x])=\{qs^{-1}(x)|x\in \operatorname{Ker} \partial\}\in
\operatorname{Ker} i^!_-/\operatorname{Im}c
$$
that is represented by an  element $q(y)$ where $y\in LP_{n-1}(F)$
and $s(y)=x$.

Let  $\Cal X$ be a filtration
$$
X_k\subset X_{k-1}\subset \cdots \subset X_2\subset X_1\subset X_0=X
\tag 2.4
$$
of a closed manifold $X^n$ by means of locally-flat closed submanifolds
such that every pair of submanifolds is a manifold pair in the sense of Ranicki \cite{15}.
A filtration in (2.4), for which every pair of submanifolds
$(X_i,X_{i+1})\ (0\leq i\leq k-1) $ is a Browder-Livesay pair, is called a
{\sl  Browder-Livesay filtration} (see \cite{1},  \cite{6},  and  \cite{13}). In what follows
we shall consider only Browder-Livesay filtrations and we shall  assume that
$\text{dim}\, \,  X_k=n-k\geq 5$.
The filtration  $\Cal X$ in  (2.4) is  a {\sl stratified manifold} in the sense of
Browder-Quinn (see \cite{4}, \cite{13}, and \cite{18}).

Let  $F_i \ (0\leq i\leq k-1)$ be
a square of fundamental groups in the splitting problem
for the manifold pair$(X_i,X_{i+1})$  of the filtration in  (2.4), $G_i=\pi_1(X_i)$, and
$\rho_i=\pi_1(X_i\setminus X_{i+1})$.
Then  $LN_*(\rho_i\to G_i)$ are the splitting obstruction groups for the manifold
pair $(X_i,X_{i+1})$.

 Every inclusion $\rho_i\to G_i$ of index 2 gives a commutative braid of exact
 sequences that is similar to (2.3).
Putting together central squares from these diagrams (see \cite{9} and \cite{11}) we obtain the
following  commutative diagram

$$
\matrix
 \overset{}\to{\longrightarrow} &   L_{n}(G_0) & \overset{\partial_0}\to{\rightarrow} &LN_{n-2}(\rho_0\to G_0) \cr
   &s\nearrow \ \ \ \ \                           \searrow  &  &\cr
  LP_{n-1}(F_0)& \Gamma\downarrow & L_{n}(\rho_0\to G_0) &\cr
    &q\searrow \ \ \ \ \   \nearrow  &  &\cr
\longrightarrow &   L_{n-1}(G_1) &
\overset{\partial_1}\to{\rightarrow}  &LN_{n-3}(\rho_1\to
G_1)\cr
   &s\nearrow \ \ \ \ \                           \searrow  &  &\cr
  LP_{n-2}(F_1)& \Gamma\downarrow & L_{n-1}(\rho_1\to G_1) &\cr
    &q\searrow \ \ \ \ \   \nearrow  &  &\cr
\longrightarrow &   L_{n-2}(G_2) &
\overset{\partial_2}\to{\rightarrow} &LN_{n}(\rho_2\to G_2) \cr
   &s\nearrow \ \ \ \ \                           \searrow  & & \cr
  LP_{n-3}(F_2)& \Gamma\downarrow & L_{n-2}(\rho_2\to G_2) &\cr
    &q\searrow \ \ \ \ \   \nearrow  & & \cr
& L_{n-3}(G_3) & &\cr
                           & \vdots&&& &&\cr
\longrightarrow &   L_{n-k+1}(G_{k-1}) &
\overset{\partial_{k-1}}\to{\rightarrow} &LN_{n-k-1}(\rho_{k-1}\to G_{k-1}) \cr
   &s\nearrow \ \ \ \ \                           \searrow  & & \cr
  LP_{n-k}(F_{k-1})& \Gamma\downarrow & L_{n-k+1}(\rho_{k-1}\to G_{k-1}) &\cr
    &q\searrow \ \ \ \ \   \nearrow  & & \cr
& L_{n-k}(G_k) & &\cr
 \endmatrix
\tag 2.5
 $$
 \smallskip

\noindent In this diagram we denote by $s$ and $q$ the  similar maps from
different diagrams. However, in what follows,  it will be clear from
the context which map is under consideration. Note that the groups
and the maps in Diagram (2.5) are defined by the subscripts taken
$\bmod \ 4$ since $L_*$@-groups and Diagram (2.3) are  four-periodic.

Now we can give an inductive definition of  the sets
$$
\Gamma^j(x)\subset L_{n-j}(G_j) \ \text{for} \ (0\leq j\leq k)
$$
and iterated Browder-Livesay $j$-invariants $(1\leq j\leq k)$ with
respect to the filtration  (2.4) (see 
\cite{6}, \cite{11},   \cite{12}, and \cite{13}).

\proclaim{Definition 2.1}  Let  $x\in L_n(G_0)$. By definition,
$$
\Gamma^0(x)=\{x\}\subset L_n(G_0).
$$
The set $\Gamma^0(x)$ said to be  trivial if $x\in \operatorname{Image\{L_n(\rho_0)\to L_n(G_0)\}}$.
 Let a set
 $$
 \Gamma^{j}(x)\subset L_{n-j}(G_{j})\ \ (0\leq j\leq k-1)
 $$
 be defined. For $j\geq 1$, it is called trivial if $0\in \Gamma^{j}(x)$.

If $\Gamma^{j}(x) (0\leq j\leq k-1)$ is defined and nontrivial, then the $(j+1)-th$
Browder-Livesay invariant with respect filtration  (2.4) is the set
$$
\partial_{j} (\Gamma^{j}(x))\subset LN_{n-j-2}(\rho_{j-1}\to
G_{j-1}).
$$
 The
$(j+1)-th$ invariant is nontrivial if $0\notin \partial_{j}
(\Gamma^{j}(x))$.

If the $(j+1)-th \ (1\leq j\leq k-1)$ Browder-Livesay invariant is
defined and trivial then the set  $\Gamma^{j+1}(x)$ is defined as
$$
\Gamma^{j+1}(x)\overset{def}\to{=}
\Gamma(\Gamma^{j}(x))\overset{def}\to{=} \{qs^{-1}(z)|z\in
\Gamma^{j}(x), \partial_{j}(z)=0\}\subset L_{n-j-1}(G_{j+1}).
$$
\endproclaim

\proclaim{Theorem 2.2} (\cite{11} and \cite{13}) Let $x\in L_n(G_0)$ be
an element with  a  nontrivial $j-th$ Browder-Livesay invariant relatively
to a Browder-Livesay filtration $\Cal X$ of the manifold $X$. Then
the element $x$ cannot be realized by a normal map of closed manifolds.
\endproclaim

Let us consider an   an infinite  diagram $\Cal D_{\infty}$ of groups
with orientations
$$
\matrix
        &         &         & \rho_2&         &   &         & \rho_1&
              &     &         &\rho_0&         &   \\
\cdots \searrow    &         &\swarrow &       &\searrow &   &\swarrow
&
     &\searrow &     &\swarrow &      &\searrow &     \\
        & G_3&         & \overset{\cong}\to{\to}      &         &G_2&
               &  \overset{\cong}\to{\to}      &         & G_1 &
                      &  \overset{\cong}\to{\to}     &         &G_0, \\
\endmatrix
\tag 2.6
$$
which is commutative as  the diagram of groups.
The maps $\rho_i\to G_i$ and $\rho_{i}\to G_{i+1}$
in (2.6) are index 2 inclusions
 of groups with orientations,
  and the horizontal maps
preserve the orientations on the images  of the groups $\rho_i$ and reverse the
orientations outside these images.  Each commutative triangle from (2.6)
defines an algebraic version of  diagram  (2.3) for the inclusion $\rho_i\to G_i$
\cite{16}. Putting together central squares of these diagrams
we obtain an infinite in bottom direction diagram which is similar to
diagram (2.5) (see \cite{1}, \cite{6}, \cite{9},  and \cite{13}).
Thus we can define the {\sl iterated Browder--Livesay
invariants  of an element
$x\in L_n(G_0)$ relatively to a diagram in (2.6)}
similarly to the case of the filtration $\Cal X$. A result 
similar to Theorem 2.2 is true for Browder-Livesay invariants
of an element $x$ relatively  to a diagram in (2.6) \cite{13}.
Note that the finite subdiagram (from $G_k$ until $G_0$)  of
diagram (2.6) provides a commutative diagram (2.5). Denote this diagram
by $\Cal D_k$.

\proclaim{Definition 2.3} (\cite{1}, \cite{6},  and \cite{12}) Let $x\in L_n(G_0)$.
The element $x$ is the element of the second type with respect to $\Cal D_{\infty}$ if
all the sets $\Gamma^j(x) \ (j\geq 0)$ are defined and nontrivial and
all Browder-Livesay invariants with respect to $\Cal D_{\infty}$ are defined and trivial.
\endproclaim

\proclaim{Theorem 2.4} (\cite{11}  and \cite{12}) Let $x\in L_{n}(G_0)$ be an element of the second type
with respect to $\Cal D_{\infty}$. Then $x$ cannot be realized by a normal map
of closed manifolds.
\endproclaim

Now we state the main results of the paper.

\proclaim{Theorem 2.5} Let $x\in L_n(G_0)$ $(n\geq 5)$ be an element
for which the set $\Gamma^3(x)$ is defined and nontrivial with
respect to the  subdiagram  $\Cal D_3$ of diagram (2.6).
Then the element $x$ does not belong to
the subgroup $I_n(\pi)$.
\endproclaim

\proclaim{Remark} To define  $j-th$ Browder-Livesay invariant of $x$  for
$j\geq 4$ relatively to diagram $\Cal D_{\infty}$ we must have a nontrivial set
$\Gamma^3(x)$ with respect to  subdiagram $\Cal D_3$ of $D_{\infty}$.
Thus,  the nontriviality of a  the $j-th$ Browder-Livesay invariant for
$j\geq 4$ automatically implies nontriviality of the   invariant $\Gamma^3(x)$.
Also, if the element $x$ has the second type, then by Definition 2,
 the   invariant $\Gamma^3(x)$ is nontrivial.
\endproclaim
\bigskip
\subhead 3. Proof of  Theorem 2.5
\endsubhead
\bigskip

\proclaim{Lemma 3.1} Let $X^n=X^n_0$ be a manifold with the fundamental group $G_0$
and $\Cal D_k$ be a finite subdiagram of diagram (2.6) such that $n-k\geq 5$.
 Then there exists a Browder-Livesay filtration as in (2.4) $\Cal X$ of $X=X_0$
which corresponds to the diagram $\Cal D_k$.
\endproclaim
\demo{Proof} Consider a  map
$$
 \phi : X^n\to RP^N
$$ to  a real projective space of high dimension which induces an
epimorphism of fundamental groups $G_0\to \Bbb Z/2$ that
has the  kernel $\rho_0$. Using the standard arguments
(see \cite{7}, \cite{9}, \cite{12}, and \cite{17, \S 11, 12C}),
 we can suppose that the map $\phi $ is
transversal to $RP^{N-1}\subset  RP^{N}$ with $ Y^{n-1}=(\phi_n)^{-1}(RP^{N-1})$
such that  the induced  map
$\pi_1(Y)\to \pi_1(X)$ is an isomorphism.
The pair $(X,Y)$ is the  Browder-Livesay pair
that gives the filtration with $\Cal D_1$. Iterating this process
we obtain the desired  result.
\enddemo

\demo{Proof of  Theorem 2.5} Let the  element $x\in L_n(G_0)$ act trivially  on a
manifold $X^{n-1}$. Taking the product with  $\Bbb RP^4$, we can suppose that
 dimension $n-1\geq 8$ (see \cite{17}). Consider a Browder-Livesay filtration
$$
X_3\subset X_2\subset X_1\subset X_0=X^{n-1}
$$
which gives the diagram $\Cal D_3$ by Lemma 1.
Let $U$ be a tubular neighborhood of $X_3$ in $X_0$. Note,
that $\pi_1(X_0\setminus X_3)=\pi_1(X_0)=G_0$.
 In accordance with  \cite{17},  we can construct the action of the element $x$ on the manifold $X_0$
"outside the tubular neighborhood $U$". The proof is identically with the proof
of Theorem 5.8 (respectively Theorem 6.5) in \cite{17} because 
$\pi_1(X\setminus U)\cong \pi_1(X)$.

 This means
that we can  represent $X_0\times I$ as
$$
X_0\times I= {(X_0\setminus U)}\times I\bigcup_{\partial U\times I} \overline{ U}\times I
$$
such that the normal map
$$
F\colon (W^{n};\partial_0W,\partial_1W)\to (X\times I;X\times\{0\},
X\times\{1\}) \tag 3.1
$$
has the following properties.

i) The manifold $W^n$ is a union
$$
W^n= V^n\bigcup_{\partial U\times I} \overline{ U}\times I
$$
where
$$
\partial V= \partial_0V\bigcup_{\partial U\times \{0\}} \partial U\times I
\bigcup_{\partial U\times \{1\}} \partial_1V,
$$
$$
\partial_0V=(X\setminus U)\times \{0\},
$$
the boundary of $ \partial_0V$
is equal to
$$
\partial  U\times \{0\},
$$
the boundary of $\partial_1V$
is equal to
$$
 \partial U\times \{1\},
$$
$$
\partial_0W = X\times\{0\}=(X\setminus U)\times \{0\}\bigcup_{\partial U\times \{0\}} \overline{ U}\times \{0\},
$$
$$
\partial_1W = \partial_1V\bigcup_{\partial U\times \{1\}} \overline{ U}\times \{1\}.
$$

ii) The restriction
$$
F|_{X\times \{0\}\bigcup \overline{ U}\times I}\colon X\times \{0\}\bigcup_{\overline{ U}\times \{0\}} \overline{ U}\times I\to
X\times \{0\}\bigcup_{\overline{ U}\times \{0\}}  \overline{ U}\times I
$$
is the identity map,  the restriction
$$
F|_{\partial_1 V}\colon \partial_1V\to (X\setminus U)\times \{1\}
$$
is a homeomorphism which  is identity map on the boundary $\partial U\times \{1\}$,
and  the restriction
$$
F|_{\partial_1 W}\colon \partial_1W\to X\times \{1\}
\tag 3.2
$$
is a homeomorphism.

\smallskip

Changing
$$
F|_V\colon  V\to (X\setminus U)\times I
$$
relatively boundary $\partial V$ in the class of normal bordisms,   we can suppose
that it is transversal to
$$
\left(X_1\setminus (X_1\cap U)\right)\times I\subset (X\setminus U)\times I
$$
with a  transversal preimage
$
(V_1, \partial V_1),
$
where
$$
\partial V_1= \partial_0V_1\bigcup_{(X_1\cap \partial U)\times \{0\}} (X_1\cap \partial U)\times I
\bigcup_{(X_1\cap \partial U)\times \{1\}}  \partial_1V_1=
$$
$$
=[X_1\setminus (X_1\cap U)]\times \{0\}\bigcup_{(X_1\cap \partial U)\times \{0\}} (X_1\cap \partial U)\times I
\bigcup_{(X_1\cap \partial U)\times \{1\}} \partial_1V_1,
$$
 the restriction of $F|_{V}$
to
$$
\partial_0V_1\bigcup_{(X_1\cap \partial U)\times \{0\}} (X_1\cap \partial U)\times I
$$
is
the identity map, and the restriction of $F|_{V}$
to $\partial_1V_1$ is a homeomorphism
$$
F|_{\partial_1V_1}\colon \partial_1V_1\to [X_1\setminus (X_1\cap U)]\times I.
$$
It follows that
 the map $F$ in (3.1)  is transversal to $X_1\times I$ with a transversal
preimage
$(W_1, \partial W_1)$ where
$\partial W_1= \partial_0W_1\cup \partial_1 W_1$. Let $U_1=U\cap X_1$ be the tubular neighborhood of
$X_3\subset X_1$. The restriction
$F_1=F|_{W_1}$ is a
normal  map
$$
F_1\colon (W^{n-1}_1;\partial_0W_1,\partial_1W_1)\to (X_1\times I;X_1\times\{0\},
X_1\times\{1\}) \tag 3.3
$$
with  the following
 properties:

i) the manifold $W^{n-1}_1$ is a union
$$
W_1= V_1\bigcup_{\partial U_1\times I} \overline{ U_1}\times I
$$
where
$$
\partial V_1= \partial_0V_1\bigcup_{\partial U_1\times \{0\}} \partial U_1\times I
\bigcup_{\partial U_1\times \{1\}} \partial_1V_1,
$$
$$
\partial_0V_1=(X_1\setminus U_1)\times \{0\},
$$
the boundary of $\partial_0V_1$
is equal to
$$
\partial  U_1\times \{0\},
$$
the boundary of $\partial_1V_1$
is equal to
$$
 \partial U_1\times \{1\},
$$
$$
\partial_0W_1 = X_1\times\{0\}=(X_1\setminus U_1)\times
\{0\}\bigcup_{\partial U_1\times \{0\}} \overline{ U_1}\times \{0\}=\partial_0V_1
\bigcup_{\partial U_1\times \{0\}} \overline{ U_1}\times \{0\},
$$
$$
\partial_1W_1 = \partial_1V_1\bigcup_{\partial U_1\times \{1\}} \overline{ U_1}\times \{1\};
$$

ii) the restriction
$$
F_1|_{X_1\times \{0\}\bigcup \overline{ U_1}\times I}\colon X_1\times \{0\}\bigcup_{\overline{ U_1}\times \{0\}}
 \overline{ U_1}\times I\to
X_1\times \{0\}\bigcup_{\overline{ U_1}\times \{0\}}  \overline{ U_1}\times I
$$
is the identity map,
the restriction
$$
F_1|_{\partial_1 V_1}\colon \partial_1V_1\to (X_1\setminus U_1)\times \{1\}
$$
is a homeomorphism which  is identity map on the boundary $\partial U_1\times \{1\}$,
and hence the restriction
$$
F_1|_{\partial_1 W_1}\colon \partial_1W_1\to X_1\times \{1\}
\tag 3.4
$$
is a homeomorphism.
\smallskip

Now consider the diagram (2.5) for which $x\in L_n(G_0)$. Since $\sigma (F)=x$ and the map in (3.2)
is a homeomorphism.   It follows from the geometrical definition of the map
$$
\partial_0:L_n(G_0)\to LN_{n-2}(\rho_0\to G_0)
$$
(see \cite{7}, \cite{12}, and \cite{17}
that  $\partial_0(x)=0$.
It follows from  geometrical definition of the map $\Gamma$ in  diagram (2.5) (\cite{11}, \cite{12} and \cite{17})
that the normal map $F_1$ in (3.3) has a surgery obstruction $\sigma (F_1)=x_1\in L_{n-1}(G_1)$
which lies in the class $\Gamma(x)$. As above, since the restriction
$
F_1|_{\partial_1 W_1}
$
in (3.4) is a homeomorphism, we obtain
$\partial_1(x_1)=0\in L_{n-3}(\rho_1\to G_1)$.
For passing from the map $F_1$ to a map $F_2$, which is restriction of $F_1$ to the a transversal
preimage of $X_2\times I$, we can use the same line of arguments as for
the passing from the map $F$ to the map  $F_1$. Let $U_2=U\cap X_2$ be the tubular neighborhood
of $X_3\subset X_2$.
Thus we obtain a normal map
$$
F_2\colon (W^{n-2}_2;\partial_0W_2,\partial_1W_2)\to (X_2\times I;X_2\times\{0\},
X_2\times\{1\}) \tag 3.5
$$
with  a decomposition
$$
W_2= V_2\bigcup_{\partial U_2\times I} \overline{ U_2}\times I
$$
which is similar to the decomposition above. In particular,
$$
\partial V_2= \partial_0V_2\bigcup_{\partial U_2\times \{0\}} \partial U_2\times I
\bigcup_{\partial U_2\times \{1\}} \partial_1V_2,
$$
$$
\partial_0V_2=(X_2\setminus U_2)\times \{0\},
$$
the boundary of $\partial_0V_2$
is equal to
$$
\partial  U_2\times \{0\},
$$
the boundary of $\partial_1V_2$
is equal to
$$
 \partial U_2\times \{1\},
$$
$$
\partial_0W_2 = X_2\times\{0\}, \ \
\partial_1W_2 = \partial_1V_2\bigcup_{\partial U_2\times \{1\}} \overline{ U_2}\times \{1\}.
$$
The map $F_2$  is the identity
on
$$
\partial_0W_2 \bigcup_{\overline{ U_2}\times \{0\}}
\overline{ U_2}\times I
=X_2\times \{0\}\bigcup_{\overline{ U_2}\times \{0\}}
 \overline{ U_2}\times I,
$$
and the  restriction of $F_2$ to $\partial_1W_2$ is a homeomorphism
$$
F_2|_{\partial_1W_2}\colon \partial_1W_2\to X_2\times \{1\}
$$
which is identity on
$$
\overline{U_2}\times\{1\}\subset {\partial_1W_2}.
$$
As above,  the normal map $F_2$ in (3.5) has a surgery obstruction $\sigma (F_2)=x_2\in L_{n-2}(G_2)$
which lies in the class $\Gamma^2(x)$ and $\partial_2(x_2)=0$.
By our construction, 
$$
F_2^{-1}(U_2\times\{I\})=U_2\times\{I\}\subset W_3
$$
and
$
F_2|_{U_2\times\{I\}}
$
is identity. Since $U_2$ is a tubular neighborhood of $X_3$ in $X_2$
we obtain that a restriction of $F_2$ to the transversal
preimage of $X_3\times I$
$$
F_3=F_2|_{F_2^{-1}(X_3\times I)}
$$
  is the  identity.
Thus, the surgery obstruction $\sigma (F_3)\in L_{n-3}(G_3)$ is trivial.
This obstruction lies in the class $\Gamma^3(x)$,  and hence
$0\in \Gamma^3(x)$ and $\Gamma^3(x)$ is trivial.
The theorem is thus  proved.
\enddemo
\smallskip

There are many examples of nontrivial first and second Browder-Livesay 
invariants and notrivial classes $\Gamma^{3}$ (see for example \cite{12}).
We do not  know examples with nontrivial third Browder-Livesay invariant. 
Also, in all known cases  nontriviality of the class 
$\Gamma^{2}(x)$ for some element $x\in L_n(\pi)$ implies that 
element $x$ does not lie in $C_n(\pi)$. It would  be very interesting 
to understand this problem.

\bigskip

Authors' addresses:

\bigskip

\noindent
Friedrich Hegenbarth: Dipartimento di Matematica, Universit\`a di
Milano,
Via Saldini n. 50, 20133 Milano, Italia;
\smallskip
\noindent
E--mail: hegenbar\@mat.unimi.it

\bigskip
\noindent
Yuri V. Muranov:  Division de Estudios de Posgrado, Universidad Tecnologica de la Mixteca,
km. 2.5 Carretera a Acatlima, Huajuapan de Leon, Oaxaca, C.P. 69000, Mexico;
\smallskip
\noindent
E--mail: ymuranov\@mixteco.utm.mx

\bigskip

\noindent Du\v san Repov\v s: Faculty of Mathematics ans  Physics, University of
Ljubljana, Jadranska 19, Ljubljana, Slovenia;
\smallskip

\noindent E--mail: dusan.repovs\@guest.arnes.si

\enddocument

\Refs

\ref\no 1
\by A. Bak -- Yu. V. Muranov
\paper Splitting a simple homotopy equivalence
  along a submanifold with filtration
\jour   Sbornik : Mathematics 
\vol 199
\issue 6 
\pages 1–23
\yr 2008 
\endref

\ref\no 2
\by A. Bak -- Yu. V. Muranov
\paper Splitting along submanifolds, and $\Bbb L$-spectra (Russian)
\jour Sovrem. Mat. Prilozh., Topol., Anal. Smezh. Vopr.
\issue 1
\yr 2003
\pages 3--18
\moreref English translation in J. Math. Sci. (N. Y.) 123 (2004), no. 4, 4169--4184
\endref

\ref\no  3
\by W. Browder -- G. R. Livesay
\paper Fixed point free involutions on
homotopy spheres
\jour Bull. Amer. Math. Soc.
\vol 73
\yr 1967
\pages 242--245
\endref

\ref\no 4
\by W. Browder -- F. Quinn
\paper A surgery theory for G-manifolds and stratified sets
\jour{\rm in}\ Manifolds--Tokyo 1973
\publ Univ. of Tokyo Press
\yr 1975
\pages 27--36
\endref


\ref\no 5
\by S. E. Cappell  -- J. L. Shaneson
\paper Pseudo-free actions. I
\jour Lecture Notes in Math.
\vol 763
\yr 1979
\pages  395--447
\endref

\ref\no 6
\by A. Cavicchioli -- Yu. V. Muranov -- F. Spaggiari
\paper On the elements ofthe  second type in surgery groups
\jour  Preprint MPI
\yr 2006
\endref


\ref\no 7
\by I. Hambleton
\paper Projective surgery obstructions
on closed manifolds
\jour Lecture Notes in Math. \yr 1982 \vol 967
\pages 101--131
\endref


\ref\no 8
\by I. Hambleton -- J. Milgram -- L. Taylor -- B. Williams
\paper Surgery with finite fundamental group
\jour Proc. London Mat. Soc.
\vol 56
\yr 1988
\pages 349--379
\endref

\ref\no 9
\by I. Hambleton -- A. F. Kharshiladze
\paper A spectral sequence in surgery theory
\jour Mat. Sbornik
\yr 1992
\vol 183
\pages 3--14
\transl\nofrills English transl. in
\jour  Russian Acad. Sci. Sb. Math.
\vol 77
\yr 1994
\pages 1--9
\endref

\ref\no 10
\by I. Hambleton -- A. Ranicki -- L. Taylor
\paper Round L-theory
\jour J.  Pure  Appl. Algebra
\vol 47
\yr 1987
\pages  131--154
\endref


\ref\no 11
\by A. F. Kharshiladze
\paper Iterated Browder-Livesay invariants and oozing problem
\jour Mat. Zametki
\vol 41
\yr 1987
\pages 557--563
\transl\nofrills English transl. in
\jour Math. Notes
\vol 41
\yr 1987
\endref

\ref\no 12
\by A. F. Kharshiladze
\paper Surgery on manifolds with finite fundamental groups
\jour Uspechi Mat. Nauk
\vol 42
\yr 1987
\pages 55--85
\transl\nofrills English transl. in
\jour Russian Math. Surveys
\vol 42
\yr 1987
\endref





\ref\no 13 \by  Yu. V. Muranov -- D. Repov\v s -- R. Jimenez
\paper Surgery spectral sequence and manifolds with filtrations
\publ Trudy MMO (in Russian) \vol 67 \yr 2006 \pages 294--325
\endref


\ref \no 14 \by A. A. Ranicki \paper The total surgery obstruction
\jour Lecture Notes in Math. \vol 763 \yr 1979 \pages 275--316
\endref


\ref \no 15 \by A. A. Ranicki \book Exact Sequences in the
Algebraic Theory of Surgery \publ Math. Notes  {\bf 26}, Princeton
Univ. Press
 \publaddr Princeton, N. J.
\yr 1981
\endref

\ref \no 16 \by A. A. Ranicki \paper The L-theory of twisted
quadratic extensions \jour Canad. J. Math. \yr 1987 \vol 39 \pages
245--364
\endref


\ref \no 17 \by C. T. C. Wall \book Surgery on Compact Manifolds
\publ Academic Press \publaddr London -- New York \yr 1970
\moreref \by \ Second Edition (A. A. Ranicki, ed.) \publ Amer.
Math. Soc., Providence, R.I. \yr 1999
\endref

\ref\no 18
 \by S. Weinberger
 \book The Topological Classification of Stratified Spaces
 \publ The University of Chicago Press
\publaddr Chicago -- London
 \yr 1994
 \endref



\endRefs
\bye